\newcommand{\rem}[1]{}
\documentclass{amsart}
\usepackage{amsfonts,amssymb,amsmath,amsthm,mathrsfs}
\usepackage{url}
\usepackage[dvips]{epsfig}
\newtheorem{thrm}{Theorem}[section]

\newtheorem{remark}[thrm]{Remark}
\theoremstyle{definition}

\begin{document}
\author[C.~A.~Mantica and L.~G.~Molinari]
{Carlo~Alberto~Mantica and Luca~Guido~Molinari}
\address{C.~A.~Mantica (corresponding author): Physics Department, Universitˆ degli Studi di Milano ,Via Celoria 16, 20133, Milano and
I.I.S. Lagrange, Via L. Modignani 65, 
20161, Milano, Italy -- L.~G.~Molinari: Physics Department,
Universit\`a degli Studi di Milano and I.N.F.N. sez. Milano,
Via Celoria 16, 20133 Milano, Italy.}
\email{carloalberto.mantica@libero.it, luca.molinari@unimi.it}
\subjclass[2010]{Primary 53B30, 53B50, Secondary 53C80, 83C15}
\keywords{extended recurrent manifold, Robertson-Walker space-time, perfect fluid space-time, 
torse-forming vector, concircular vector.}
\title[Extended recurrent Lorentzian manifolds]
{A note on extended recurrent\\ Lorentzian manifolds }
\begin{abstract} 
Extended recurrent pseudo-Riemannian manifolds were introduced by Mileva Prvanovi\'c. 
We reconsider her work in the light of recent results and show that the manifold is conformally flat, and it is a space of quasi-constant curvature. 
We also show that an extended recurrent Lorentzian manifold, with time-like associated covector, is a 
perfect fluid Robertson-Walker space-time. We obtain the equation of state; in $n=4$ and if the scalar curvature 
is zero, a model for incoherent radiation is obtained.
\end{abstract}
\date{25 aug 2016}
\maketitle
\centerline{\sf Dedicated to the memory of Dr. Mileva Prvanovi\'c}
{\quad}\\
\section{Introduction}
In 1999 Mileva Prvanovi\'c \cite{[20]} introduced the following differential structure on a pseudo-Riemannian manifold, that she
named ``extended recurrent manifold'':
\begin{align} \label{eq_1.1}
& \nabla_i R_{jklm} = A_i R_{jklm} +(\beta -\psi) A_iG_{jklm} \\
&+ \frac{\beta}{2} [A_j G_{iklm}+A_k G_{jilm}+A_l G_{jkim}+A_m G_{jkli}] \nonumber
\end{align}
$A_i$ is a closed one-form named ``associated covector'', $\beta $ and $\psi$ are non vanishing scalar functions with $\nabla_j\psi = \beta A_j$,  
 $G_{jklm}= g_{mj} g_{kl} - g_{mk} g_{jl}$.\\ 
She proved that the associated covector is a concircular vector: $\nabla_s A_r = f g_{rs} +hA_r A_s$ with scalar functions $f$ and $h$, and showed that the metric has the warped form
\begin{align}\label{eq_1.2}
ds^2 = (dx^1)^2 + e^{\eta} g^*_{\alpha\beta} dx^\alpha dx^\beta
\end{align}
where $g^*_{\alpha\beta}$ are functions only of $x^\gamma $ ($\gamma= 2,\dots, n$) and $\eta $ is a scalar function of $x^1$.\\
These properties will be reviewed in Section 2, where we also derive some new ones.
In particular we show that an
extended recurrent pseudo-Riemannian manifold is conformally flat, and it is a space of quasi
constant curvature, according the definition by K. Yano and B.-Y. Chen \cite{[5]}. 
In Section 3 we focus on extended recurrent Lorentzian manifolds (space-times).
Based on our recent study of Generalized Robertson Walker manifolds, to which the present
model eventually belongs, we show that an extended recurrent space-time with
time-like associated covector is a perfect fluid Robertson-Walker spacetime.  The barotropic equation of state is obtained; in the particular case of vanishing scalar curvature, 
in 4 dimensions, we obtain a model for incoherent radiation. 

Throughout the paper we adopt the convention $R_{ij} = R_{imj}{}^m$ and $R=R^m{}_m$ for the Ricci tensor and the scalar curvature, and use the notation $v^2 = v^mv_m$.

\section{General properties of extended recurrent\\ pseudo-Riemannian manifolds}
We review some basic properties of extended recurrent pseudo-Riemannian
manifolds exposed in \cite{[20]}. Furthermore, we prove some new characterizations of such manifolds.

Following the procedure in \cite{[20]}, by contracting \eqref{eq_1.1} with $g^{jm}$ it is
\begin{align}\label{eq_2.1}
\nabla_i R_{kl} = A_i [R_{kl} -g_{kl}(n\beta -(n-1)\psi ) ] -\frac{\beta}{2}(n-2)(A_k g_{il} +A_l g_{ik}).
\end{align}
Contracting again \eqref{eq_2.1} with $g^{kl}$  we obtain
\begin{align}\label{eq_2.2}
\nabla_i R= A_i [R-(n^2+n-2)\beta + n(n-1)\psi].
\end{align}
On the other hand, by the second Bianchi identity for the Riemann tensor it is
$A_i (R_{jklm}-\psi G_{jklm}) +A_j (R_{kilm}-\psi G_{kilm}) +A_k (R_{ijlm}-\psi G_{ijlm})=0 $. Contracting this with $g^{im}$ it is
\begin{align}\label{eq_2.4}
R_{jklm} A^m= A_k [R_{jl}+\psi (n-2)g_{jl}] -A_j [R_{kl}+\psi (n-2) g_{kl}].
\end{align}
and contracting \eqref{eq_2.4} with $g^{kl}$ we obtain
\begin{align}\label{eq_2.5}
R_{jm} A^m = \tfrac{1}{2}A_j [R+\psi (n-2)(n-1)].
\end{align}
The components of the Weyl conformal curvature tensor are \cite{[19]}:
\begin{align}\label{eq_2.3}
C_{jkl}{}^m = &R_{jkl}{}^m +\frac{1}{n-2}(g_{jm}R_{kl}-g_{km}R_{jl} + R_{jm}g_{kl}-R_{km}g_{jl} ) \\
&-\frac{g_{jm}g_{kl}-g_{km}g_{jl} }{(n-1)(n-2)} R\nonumber 
\end{align}
By taking the covariant derivative of \eqref{eq_2.3} 
and inserting \eqref{eq_2.2} and \eqref{eq_2.1} we infer that 
\begin{align}\label{nablaC=AC}
\nabla_i C_{jklm}= A_i C_{jklm}
\end{align} 
Now, \eqref{eq_2.4}, \eqref{eq_2.5} are used to evaluate $A_m C_{jkl}{}^m$:
\begin{align}\label{eq_2.6}
 A_m C_{jkl}{}^m = \frac{n-3}{n-2} \Big [& A_k \left ( R_{jl} -\frac{R-\psi (n-1)(n-2)}{2(n-1)} g_{jl}\right )\\
 &-A_j \left (R_{kl} - \frac{R-\psi (n-1)(n-2)}{2(n-1)}g_{kl}\right )\Big ]  \nonumber
\end{align}
Next, consider Lovelock's identity (\cite{[14]} page 289):
\begin{align*} 
\nabla_i \nabla_m R_{jkl}{}^m + \nabla_j \nabla_m R_{kil}{}^m +\nabla_k \nabla_m R_{ijl}{}^m = -R_{im}R_{jkl}{}^m - R_{jm}R_{kil}{}^m - R_{km} R_{ijl}{}^m
\end{align*}
The evaluation of  $\nabla_i \nabla_m R_{jkl}{}^m + \nabla_j \nabla_m R_{kil}{}^m +\nabla_k \nabla_m R_{ijl}{}^m$ with the aid of \eqref{eq_2.1} 
gives zero, therefore it is $R_{im}R_{jkl}{}^m + R_{jm}R_{kil}{}^m + R_{km} R_{ijl}{}^m =0$.
By taking the covariant derivative $\nabla_s$ of the last
expression and contracting with $g^{is}$, after long calculations, it is inferred that (provided $\beta\neq 0$ and $n>3$)
\begin{align}\label{eq_2.8} 
A_j \left [ R_{kl} - g_{kl} \frac{R-\psi(n-1)(n-2)}{2(n-1)} \right ] = A_k \left [ R_{jl} - g_{jl} \frac{R-\psi (n-1)(n-2)}{2(n-1)} \right ] 
\end{align}
From \eqref{eq_2.8} and \eqref{eq_2.6} immediately it is 
$\nabla_m C_{jkl}{}^m = A_m C_{jkl}{}^m =0$.\\
The second Bianchi identity for the Weyl tensor is (see \cite{[1]})
\begin{align*}
\nabla_i C_{jkl}{}^m +\nabla_j C_{kil}{}^m + \nabla_k C_{ijl}{}^m = \frac{1}{n-3} \Big [ \delta_j^m \nabla_p C_{kil}{}^p + \delta_k^m \nabla_p C_{ijl}{}^p \\
+ \delta_i^m \nabla_p C_{jkl}{}^p + g_{kl} \nabla_p C_{ji}{}^{mp}+ 
g_{il} \nabla_p C_{kj}{}^{mp}+ g_{jl} \nabla_p C_{ik}{}^{mp} \Big] \nonumber
\end{align*}
For a conformally recurrent manifold it becomes 
\begin{align}
A_i C_{jklm} + A_j C_{kilm}+A_k C_{ijlm} = \frac{A^p}{n-3}\Big [g_{mj} C_{kilp} + g_{mk} C_{ijlp}\\
+ g_{mj}C_{jklp}+ g_{kl}C_{jimp}+ g_{il} C_{kjmp}+g_{jl}C_{ikmp}\Big ]=0 \nonumber
\end{align}
because $A_p C_{jkl}{}^p=0$.
Thus in our case it is $A_i C_{jklm} + A_j C_{kilm}+A_k C_{ijlm} =0$ from which $A^2 C_{jklm}=0$.
Therefore, if $A^2\neq 0$, the manifold is conformally flat: $C_{jklm}=0$. Moreover if $A^2\neq 0$ eq.\eqref{eq_2.8} readily rewrites as:
\begin{align} \label{eq_2.10}
2(n-1) R_{kl} -g_{kl}(R-\psi (n-1)(n-2)) = \frac{A_kA_l}{A^2} (n-2) [ R+\psi n(n-1)]
\end{align}
and shows that the space is quasi-Einstein (see for example \cite{[8],[10],[11],[12],Prvan90}):
\begin{align}\label{Ricci}
R_{kl}=a g_{kl}+b \frac{A_kA_l}{A^2}, \quad a=\frac{R-\psi (n-1)(n-2)}{2(n-1)}, 
\; b=\frac{n-2}{2(n-1)}[R+\psi n (n-1)]
\end{align}
Inserting this in \eqref{eq_2.3} with $C_{jklm}=0$ gives the Riemann tensor:
\begin{align}\label{eq_2.11}
R_{jklm} =&\frac{b}{n-2}\left [-g_{jm} \frac{A_k A_l}{A^2} + g_{km} \frac{A_j A_l}{A^2} -g_{kl}\frac{A_j A_m}{A^2}+g_{jl}\frac{A_k A_m}{A^2}\right] \\
&+\psi (g_{jm}g_{kl}-g_{jl}g_{km}). \nonumber
\end{align}
Eq.\eqref{eq_2.11} characterizes the ``manifolds of quasi constant curvature", introduced by Chen and Yano in 1972 \cite{[5]}.
We thus proved the following
\begin{thrm}\label{thrm_1}
An $n\ge 3$ dimensional extended recurrent pseudo-Riemannian manifold is conformally flat and is 
a space of quasi-constant curvature.
\end{thrm}
Note that the hypothesis $\nabla_j \psi = A_j \beta $ is not used in the proof of Theorem 
\ref{thrm_1}.\\
As shown in \cite{[20]}, the covariant derivative $\nabla_s$ of \eqref{eq_2.10} and the condition 
$\nabla_j \psi = A_j  \beta$ imply that 
$$\nabla_s A_r = fg_{rs} + \omega_s A_r $$ 
where $f=-\frac{(n-1)\beta}{R+n(n-1)\psi} A^2$, $\omega_s =hA_s$, 
$h=\frac{A^jA^l\nabla_jA_l}{A^4}+\frac{(n-1)\beta}{R+n(n-1)\psi} $.  By showing 
$\nabla_s h =\mu A_s$ it follows that $\omega_s $ is closed (i.e. $A_j$ is a proper concircular vector). Based on the works \cite{[27],[28]} by Yano, Prvanovic in \cite{[20]}
concluded that the metric has the warped form \eqref{eq_1.2}.
\section{Extended recurrent space-times}
In this section we consider extended recurrent Lorentzian manifolds (i.e. space-times) with a time-like
associated covector ($A^2<0$). We prove it that it is a Robertson-Walker space-time. For this, 
we need a generalization of such spaces: \\
An $n\ge 3$ dimensional Lorentzian manifold is named generalized Robertson-Walker space-time (for short GRW) 
if the metric may take the shape:
\begin{align}\label{eq_3.1} 
ds^2 = -(dx^1)^2 + q(x^1)^2 g^*_{\alpha\beta} (x^2,\dots, x^n) dx^\alpha dx^\beta,
\end{align}
A GRW space-time is thus the warped product $1\times q^2 M^*$
(\cite{[2],[3],[23],[24]}) where $M^*$ is a $(n-1)-$dimensional Riemannian manifold. If $M^*$ is a 3-dimensional Riemannian
manifold of constant curvature, the space-time is called Robertson-Walker space-time. GRW space-times
are thus a wide generalization of Robertson-Walker space-times on which standard cosmology is modelled and include the Einstein-de Sitter space-time, the Friedman cosmological models, 
the static Einstein space-time, and the de Sitter space-time. 
They are inhomogeneous space-times admitting an isotropic radiation (see S\'anchez \cite{[23]}). 
We refer to the works by Romero et al. \cite{[21],[22]},
S\'anchez \cite{[23]} and Guti\'errez and Olea \cite{[13]} for an exhaustive presentation of geometric and physical properties.\\
Recently, perfect fluids with the condition $\nabla_m C_{jkl}{}^m=0$ were studied in \cite{[15]} and \cite{[16]}, where the authors showed that such spaces are GRW space-times.\\
The following deep result was proved by Bang Yen Chen, in ref.\cite{[4]} 
(for similar results see also the works by Yano \cite{[27],[28]}, Prvanovi\'c \cite{Prvan95}, and the recent paper \cite{[9]}).
\begin{thrm}[Chen]\label{thrm_2}
Let $(M,g)$ an $n\ge 3$ dimensional Lorentzian manifold. The space-time is a generalized Robertson-Walker space-time 
if and only if it admits a time-like vector of the form $\nabla_k X_j = \rho g_{kj}$.
\end{thrm}

In the previous section we reviewed Prvanovic's result that the associate covector is concircular,
$\nabla_j A_k = fg_{jk}+\omega_j A_k $, with $\omega_j =h A_j$ being a closed one-form. In this case $\omega_j =\nabla_j \sigma $ for a suitable scalar function.\\
If the associated covector is time-like, i.e. $A^2<0$ (with Lorentzian signature), then 
it can be rescaled to a time-like vector $X_k =A_k e^{-\sigma}$ such that $\nabla_j X_k = \rho g_{kj}$: in fact it is $\nabla_j X_k =(\nabla_j A_k -\omega_j A_k)e^{-\sigma} = (fe^{-\sigma}) X_k$.
By Chen's theorem \ref{thrm_2}, the space is a GRW space-time (see \cite{[15],[16]}).\\ 
Thus for $A^2<0$ Prvanovi\'c's model \eqref{eq_1.1}  is a quasi-Einstein GRW space-time with $C_{jklm}=0$. It is well known (see \cite{[7]}) that in this case the fiber is a space of
constant curvature and the GRW space-time reduces to an ordinary Robertson-Walker model. Moreover in
the region $A^2<0$, on defining $u_k=A_k/\sqrt{-A^2}$, it is $u^2=-1$ and the Ricci tensor 
\eqref{Ricci} becomes $R_{kl}=a g_{kl} - b u_k u_l$.  With this form of the Ricci tensor, a Lorentzian manifold is named perfect fluid space-time \cite{[15]}.
\begin{thrm}
An $n>3$ dimensional extended recurrent Lorentzian manifold with $A^2<0$ is a Robertson-Walker
space-time.
\end{thrm}
\begin{remark}
In \cite{[17]}, we proved that for a GRW space-time the condition $\nabla_m C_{jkl}{}^m=0$ is equivalent to have $R_{kl}=ag_{kl}+b\frac{X_kX_l}{X^2}$ where $X_j$ is the concircular vector of Chen's theorem. Prvanovi\'c's model matches these conditions.
\end{remark}

Some physical consequences are now outlined. Let $(M, g)$ be an $n$-dimensional Lorentzian manifold equipped with Einstein's field equations without cosmological constant,
\begin{align}\label{eq_3.2} 
R_{kl} -\tfrac{1}{2}R g_{kl} = \kappa T_{kl}
\end{align}
$\kappa =8\pi G$ is Einstein's gravitational constant (in units $c = 1$) and $T_{kl}$ is the stress-energy tensor describing the matter content of the space-time (see for example \cite{[6],[18],[26]}). Eq.\eqref{eq_3.2} is used to evaluate $T_{kl}$ obtaining:
\begin{align*} 
\kappa T_{kl} = -\frac{n-2}{2(n-1)} [ R+\psi (n-1)] \left ( g_{kl} - u_ku_l \right )
\end{align*}
We recognize a perfect fluid stress-energy tensor $T_{kl} =  (p+\mu)u_ku_l + p g_{kl}$, 
being $p$ the isotropic pressure, $\mu $ the energy density and $u_j$ the fluid flow velocity. It is
\begin{align}
\kappa p = -\frac{n-2}{2(n-1)}[R+\psi (n-1)], \quad \kappa \mu = - \frac{1}{2}\psi (n-1)(n-2)
\end{align}
One reads that the (non constant) function $\psi $
controls the energy density of the perfect fluid (then it must be negative). An equation of state can be written:
$$ p= \frac{\mu}{n-1} - \frac{n-2}{2(n-1)} \frac{R}{\kappa}. $$
In $n=4$ dimensions with the particular choice $R=0$, we have a model for incoherent radiation:
$p=\mu/3$ \cite{[25]} (a superposition of waves of a massless field with random propagation directions). 
%


\begin{thebibliography}{99}
%
\bibitem{[1]}
T.~Adati, and T.~Miyazawa, {\em On a Riemannian space with recurrent conformal curvature}, Tensor (N.S.) {\bf 18}
(1967), 348--354.
%
\bibitem{[2]}
 L.~Al\'{\i}as, A.~Romero, and M.~S\'anchez, {\em 
 Uniqueness of complete spacelike hypersurfaces of constant mean
curvature in generalized Robertson-Walker space-times}, Gen. Relat. Gravit. {\bf 27} n.1 (1995), 71--84.
%
\bibitem{[3]}
 L.~Al\'as, A.~Romero, M.~S\'anchez, Compact spacelike hypersurfaces of constant mean curvature in generalized Robertson-Walker spacetimes. In: Dillen F. editor. Geometry and Topology of Submanifolds VII. River Edge NJ, USA: World Scientific, 1995, pp 67--70.
 %
\bibitem{[4]}
 B.-Y.~Chen, {\em A simple characterization of generalized Robertson-Walker space-times}, 
 Gen. Rel. Grav. {\bf 46} (2014), 1833.
%
\bibitem{[5]} 
B.-Y.~Chen and K.~Yano, {\em Hypersurfaces of conformally flat spaces}, Tensor (N.S.) {\bf 26} (1972), 318--322.
%
\bibitem{[6]} 
J.~K.~Beem, P.~E.~Ehrlich, and K.~L.~Easley, Global Lorentzian Geometry, 2nd ed. Pure and Applied mathematics, vol. 202, 1996 Marcel Dekker, New York.
%
\bibitem{[7]}
M.~Brozos-V\'azquez, E.~Garcia-Rio, and R.~V\'azquez-Lorenzo, {\em Some remarks on locally conformally flat static
space-times}, J. Math. Phys. {\bf 46} (2005), 022501.
%
\bibitem{[8]}
M.~C.~Chaki and R.~K.~Maity, {\em On quasi-Einstein manifolds}, Publ. Math. Debrecen {\bf 57} (2000), 257--306.
%
\bibitem{[9]} 
A.~De, C.~Ozg\"ur, U.~C.~De, {\em On conformally flat Pseudo-Ricci Symmetric Spacetimes}, Int. J. Theor. Phys. {\bf 51} n.9 (2012), 2878--2887.
%
\bibitem{[10]} 
R.~Deszcz, F.~Dillen, L.~Verstraelen and L.~Vrancken, {\em Quasi-Einstein totally real submanifolds of the nearly K\"ahler 6-sphere}, 
Tohoku Math. J. {\bf 51} n.4 (1999), 461--478.
%
\bibitem{[11]} 
R.~Deszcz, M.~G{\l}ogowska, M.~Hotlo\'s and Z.~Sent\"urk, {\em On certain quasi-Einstein semisymmetric hypersurfaces},
Annu. Univ. Sci. Budapest E{\"o}tv{\"o}s Sect. Math. {\bf 41} (1998), 151--164.
%
\bibitem{[12]} 
R.~Deszcz, M.~Hotlo\'s and Z.~Sent\"urk, {\em Quasi-Einstein hypersurfaces in semi-Riemannian space forms}, Colloq. Math. {\bf  89} n.1 (2001), 81--97.
%
\bibitem{[13]} 
M.~Guti\'errez and B.~Olea, {\em Global decomposition of a Lorentzian manifold as a generalized Robertson-Walker space}, Differ. Geom. Appl. {\bf 27} (2009), 146--156.
%
\bibitem{[14]}
D.~Lovelock and H.~Rund, Tensors, Differential Forms and Variational Principles, Reprinted Edition (Dover, 1988).
%
\bibitem{[15]}
C.~A.~Mantica, L.~G.~Molinari and U.~C.~De, 
{\em A condition for a perfect fluid space-time to be a generalized Robertson-Walker space-time},
J. Math Phys. {\bf 57} n.2 (2016), 022508, Erratum,  J.M.P. {\bf 57} (2016) 049901.
%
\bibitem{[16]}
C.~A.~Mantica, Y.~J.~Suh, and U.~C.~De, 
{\em A note on generalized Robertson-Walker space-times}, Int. J. Geom. Meth. Mod. Phys. {\bf 13} (2016), 1650079, (9 pp).
%
\bibitem{[17]}
C.~A.~Mantica and L.~G.~Molinari, {\em On the Weyl and the Ricci tensors of Generalized Robertson-Walker space-times}, arXiv:1608.01209v1 [math.ph] 3 Aug 2016.
%
\bibitem{[18]}
B.~O'Neil, Semi Riemannian Geometry with applications to the Relativity (Academic Press, New York, 1983).
%
\bibitem{[19]}
M.~M.~Postnikov, Geometry VI, Riemannian geometry, Encyclopaedia of Mathematical Sciences, Vol. 91, 2001, Springer-Verlag, Berlin.
(translated from the 1998 Russian edition by S.A. Vakhrameev).
%
\bibitem{Prvan90}
M.~Prvanovi\'c, {\em On a class of SP-Sasakian manifold}, Note di Matematica, {\bf 10} n.2 (1990), 325--334.
%
\bibitem{Prvan95}
M.~Prvanovi\'c, {\em On warped product manifolds}, Filomat (Ni\v{s}) {\bf 9} n.2 (1995), 169--185.
%
\bibitem{[20]}
M.~Prvanovi\'c, {\em Extended recurrent manifolds}, Izv. Vyssh. Uchebn. Zaved. Mat. n1 (440) (1999), 41--50.
%
\bibitem{[21]} 
A.~Romero, R.~N.~Rubio, and J.~J.~Salamanca, {\em Uniqueness of complete maximal hypersurfaces in spatially parabolic generalized Robertson-Walker 
space-times}, Class. Quantum Grav. {\bf 30} n.11 (2013), 115007.
%
\bibitem{[22]} 
A.~Romero, R.~N.~Rubio, and J.~J.~Salamanca, {\em Uniqueness of complete maximal hypersurfaces in spatially parabolic generalized Robertson-Walker space-times. Applications to uniqueness results}, Int. J. Geom. Meth. Mod. Phys. {\bf 10} n.8 (2013),1360014.
%
\bibitem{[23]} 
M.~S\'anchez, {\em On the geometry of generalized Robertson-Walker spacetimes: geodesics}, Gen. Relativ. Grav. {\bf 30} (1998), 915--932.
%
\bibitem{[24]} 
M.~S\'anchez, {\em On the geometry of generalized Robertson-Walker spacetimes: curvature and Killing fields}, Gen. Relativ. Grav. {\bf 31} (1999), 1--15.
%
\bibitem{[25]} 
H.~Sthepani, D.~Kramer, M.~MacCallum, C.~Hoenselaers and E.~Hertl, Exact solutions of Einstein's Field Equations, Cambridge Monographs on Mathematical Physics 2nd ed. (Cambridge University Press, 2003).
%
\bibitem{[26]} 
R.~M.~Wald, General Relativity, (The University of Chicago Press, 1984).
%
\bibitem{[27]} 
K.~Yano, {\em Concircular geometry I-IV}, Proc. Imp. Acad. Tokyo {\bf 16} (1940), 195--200, 354--360, 442--448, 505--511.
%
\bibitem{[28]} 
K.~Yano, {\em On the torseforming direction in Riemannian Spaces}, Proc. Imp. Acad. Tokyo {\bf 20} (1944), 340--345.
\end{thebibliography}
\end{document}